\documentclass[12pt]{article}
\usepackage{amsfonts, amssymb, amsmath}

\usepackage{enumerate} 

\usepackage{scalerel}

\usepackage[colorlinks=true, urlcolor=blue,linkcolor=blue, citecolor=blue]{hyperref}

\usepackage{graphics,graphicx,amsmath}

\newcommand{\R}{\mathbb{R}}

\newcommand{\E}{\mathbb{E}}

\newcommand{\N}{\mathbb{N}}

\newtheorem{prop}{Proposition}[section]
\newtheorem{lemma}[prop]{Lemma}

\newtheorem{corollary}[prop]{Corollary}
\newtheorem{theorem}[prop]{Theorem}

\def\ind{{\rm 1\hspace{-0.90ex}1}}

\def\({\left(}
\def\){\right)}

\def\[{\left[}
\def\]{\right]}
\def\real{{\mathord{\mathbb R}}}
\def\inte{{\mathord{\mathbb N}}}
\def\Dom{\mathrm{Dom}}
\def\Var{\mathrm{Var}}

\newenvironment{Proof}{\removelastskip\par\medskip
\noindent{\em Proof.} \rm}{\penalty-20\null\hfill$\square$\par\medbreak}

\allowdisplaybreaks

\numberwithin{equation}{section}

\begin{document}
\title{
\huge
Normal approximation for sums of discrete $U$-statistics -
application to Kolmogorov bounds in random subgraph counting 
} 

\author{Nicolas Privault\thanks{
    School of Physical and Mathematical Sciences,
    Nanyang Technological University, SPMS-MAS-05-43, 21 Nanyang Link, 
Singapore 637371. e-mail: {\tt nprivault@ntu.edu.sg}. 
}
\and 
Grzegorz Serafin\thanks{Faculty of Pure and Applied Mathematics, Wroc{\l}aw University of Science and Technology, Ul. Wybrze\.ze Wyspia\'nskiego 27, Wroc{\l}aw, Poland.
e-mail: {\tt grzegorz.serafin@pwr.edu.pl}.}}

\maketitle

\vspace{-0.4cm}

\begin{abstract} 
  We derive normal approximation bounds in the Kolmogorov distance
  for sums of discrete multiple integrals and $U$-statistics
  made of independent Bernoulli random variables. 
  Such bounds are applied to normal approximation for
  the renormalized subgraphs counts in the Erd{\H o}s-R\'enyi random graph.
  This approach completely solves a long-standing conjecture in the general
  setting of arbitrary graph counting,   
  while recovering and improving recent results derived for triangles  
  as well as results using the Wasserstein distance. 
\end{abstract}

\baselineskip0.7cm

\section{Introduction}
The Mallavin approach to the Stein method
for discrete Bernoulli sequences has recently
been developed in
\cite{nourdin3}, \cite{reichenbachs}, \cite{krokowskicosa},
\cite{privaulttorrisi4}, \cite{reichenbachsAoP}, 
as an extension of the Malliavin approach to the Stein method
introduced in \cite{nourdinpeccati} for Gaussian fields. 
\\

In this paper we develop the use of multiple stochastic
integral expansions for the 
derivation of bounds on the distances between probability laws
by the Malliavin approach to the Stein and Stein-Chen methods. 
Using results of \cite{reichenbachsAoP} for general functionals
of discrete i.i.d. renormalized Bernoulli sequences
$(Y_n)_{n\in \inte}$,
we derive a Kolmogorov distance bound to the normal distribution
for sums of $U$-statistics (or multiple stochastic integrals)
of the form 
$$
\sum_{k=1}^n
\sum_{i_1,\ldots,i_k \in\inte 
  \atop
  i_r\not= i_s, \ \! \! 1\leq r\not= s \leq k}
f_k(i_1,\ldots,i_k)Y_{i_1}\cdots Y_{i_k}, 
$$
where $(Y_k)_{k\in \inte}$ is a normalized sequence of
Bernoulli random variables, 
see Theorem~\ref{prop:dKsumI}. 
We note that 
on the Erd{\H o}s-R\'enyi random graph $\mathbb{G}_n(p_n)$ 
constructed by independently retaining any edge in 
the complete graph $K_n$ on $n$ vertices
with probability $p_n\in(0,1)$, 
various random functionals admit such representations as sums of multiple integrals.
This includes the number of vertices of a given degree,
and the count of subgraphs that are isomorphic to an arbitrary
graph. 
\\
 
Our second goal is to apply such results to
the normal approximation of the renormalized
count of the subgraphs in $\mathbb{G}_n(p_n)$
which are isomorphic to an arbitrary
graph. 
Necessary and sufficient conditions for the asymptotic normality
 of the renormalization 
$$\widetilde{N}^G_n:=\frac{N^G_n-\E[N^G_n]}{\sqrt{\Var[N^G_n]}}, 
$$
 where $N^G_n$ is the number of graphs in $\mathbb{G}_n(p_n)$
 that are isomorphic to a fixed graph $G$,
 have been obtained in \cite{rucinski} where it is shown that  
\begin{align}\label{eq:conv}
  \widetilde{N}^G_n\stackrel{\mathcal{D}}{\longrightarrow}\mathcal{N} \ \mbox{ iff }
  \ \ np_n^\beta\rightarrow \infty \ \mbox{ and }\ n^2(1-p_n)\rightarrow\infty,
\end{align}
as $n$ tends to infinity,
where $\mathcal{N}$ denotes the standard normal distribution,
$$
\beta :=\max \{e_H/v_H \ : \ H\subset G\},
$$ 
and $e_H$, $v_H$ respectively denote
the numbers of edges and vertices in the graph $H$.  
\\
 
Those results have been made more precise in \cite{BKR} by
the derivation of explicit convergence rates in the Wasserstein distance
$$
 d_W (F,G):
 =\sup_{h\in\mathrm{Lip}(1)}|\mathrm{E}[h(F)]-\mathrm{E}[h(G)]|,
$$
 between the laws of random variables $F$, $G$,
 where $\mathrm{Lip}(1)$ denotes the class of real-valued
 Lipschitz functions with
 Lipschitz constant less than or equal to $1$.
In the particular case where the graph $G$ is a triangle, 
such bounds have been recently strengthened in 
\cite{roellin2} using the Kolmogorov distance
$$
 d_K (F,G):
 = \sup_{x\in \real}|P(F\leq x) - P(G\leq x)|, 
$$
 which satisfies the bound
 $d_K (F,\mathcal{N} ) \leq \sqrt{ d_W (F,\mathcal{N} )}$.
 Still in the case of triangles, 
Kolmogorov distance bounds had also been obtained
by the Malliavin approach to
the Stein method for discrete Bernoulli sequences
in
\cite{reichenbachsAoP} 
when $p_n$ takes the form $p_n = n^{-\alpha}$,
$\alpha \in [0,1)$. 
\\

In this paper we refine the results of \cite{BKR}
by using the Kolmogorov distance instead of the
Wasserstein distance.
As in \cite{BKR} we are able to consider any graph 
$G$, and therefore our results extend those
of both \cite{reichenbachsAoP} and \cite{roellin2}
which only cover the case where $G$ is a triangle.
Instead of using second order Poincar\'e inequalities
\cite{lastpeccatipenrose}, 
our method relies on an application of
Proposition~4.1 in \cite{reichenbachs}
to derive Stein approximation bounds 
for sums of multiple stochastic integrals.
\\

Our second main result Theorem~\ref{thm:main}
is a bound for the Kolmogorov 
distance between the normal distribution and the
renormalized graph count $\widetilde{N}^G_n$. 
Namely, we show that when $G$ is a graph without isolated vertices
it holds that 
\begin{equation} 
  \label{eq:main}
  d_K (\tilde{N}_G,\mathcal{N} )
  \leq C_G
  \((1-p_n)\min_{\substack{ H\subset G\\e_H\geq1}}\{n^{v_H}p_n^{e_H}\}\)^{-1/2}, 
\end{equation} 
see Theorem~\ref{thm:main}, where
$C_G >0$ is a constant depending only on $e_G$,
which improves on the Wasserstein estimates of \cite{BKR}, 
see Theorem~2 therein. 
This result relies on the representation of combined subgraph counts as
finite sums of multiple stochastic integrals, 
see Lemma~\ref{lemma4.1}, together with the application
of Theorem~\ref{prop:dKsumI} on Kolmogorov distance bounds. 
In the sequel, given two positive sequences
$(x_n)_{n\in\inte}$ and $(y_n)_{n\in\inte}$ 
we write $x_n \approx y_n$ whenever $c_1 <  x_n/y_n < c_2$
for some $c_1,c_2>0$ and all $n\in \inte$,
and for $f$ and $g$ two positive functions
we also write $f\lesssim g$ whenever $f \leq C_G g$ for some 
constant $C_G>0$ depending only on $G$. 
\\ 

Using the equivalence 
\begin{equation}
  \label{varng} 
    \Var\big[N^G_n\big]\approx
   (1-p_n)
    \max_{\substack{H\subset G \\ e_H\geq1}} \{n^{2v_G-v_H}p_n^{2e_G-e_H}\}
\end{equation} 
    as $n$ tends to infinity, see
     Lemma 3.5 in \cite{JLR}, 
     the bound \eqref{eq:main} can be rewritten in terms of the 
     variance $\Var\big[N^G_n\big]$ as 
\begin{equation} 
  \label{ddww}
  d_K \big(\widetilde{N}^G_n,\mathcal{N} \big)\lesssim
  \frac{\sqrt{\Var\big[N_G\big]}}{(1-p_n)n^{v_G}p_n^{e_G}}. 
\end{equation} 
Note that when $p_n$ is bounded away from $0$, the bound 
\eqref{eq:main} takes the simpler form
\begin{equation}
  \label{ddk} 
d_K \big(\widetilde{N}^G_n,\mathcal{N} \big)
\lesssim\frac1{n\sqrt{1-p_n}}. 
\end{equation} 
In Corollaries~\ref{cor1.3}, \ref{cor1.4} and \ref{cor1.5}
we deal with examples of subgraphs such as
cycle graphs and complete graphs, which include triangles
as particular cases, and trees. 
\\
 
In the particular case where the graph $G$ is a triangle, 
the next consequence of \eqref{eq:main} and \eqref{ddk}  
recovers the main result of \cite{roellin2},
see Theorem~1.1 therein. 
\begin{corollary}
\label{c}
For any $c\in (0,1)$, the normalized number $\widetilde{N}^G_n$ of the subgraphs in
  $\mathbb{G}_n(p_n)$ that are isomorphic to a triangle 
  satisfies 
\begin{align*}
  d_K \big(\widetilde{N}^G_n,\mathcal{N} \big)\lesssim
    \left\{\begin{array}{ll}
  \displaystyle
  \frac{1}{n\sqrt{1-p_n}} &
  \displaystyle
  \mbox{if } \ c < p_n<1,\\
  \\
  \displaystyle
 \frac{1}{n\sqrt{p_n}} &
    \displaystyle
    \mbox{if } \ n^{-1/2} < p_n \leq c, 
    \\
    \\
\displaystyle
\frac{1}{( np_n)^{3/2}} &
\displaystyle
\mbox{if } \  0 < p_n \leq n^{-1/2}.
    \end{array}\right.
\end{align*}
\end{corollary}
When $p_n$ takes the form $p_n = n^{-\alpha}$,
$\alpha \in [0,1)$, Corollary~\ref{c} 
similarly improves on the convergence rates
obtained in Theorem~1.1 of \cite{reichenbachsAoP}.
\\

\noindent 
 This paper is organized as follows. 
 In Section~\ref{s1} we recall the construction
 of random functionals of Bernoulli variables,
 together with the construction of the associated
 finite difference operator 
 and their application to Kolmogorov distance bounds
 obtained in \cite{reichenbachs}.
 In Section~\ref{s2} we derive general Kolmogorov distance 
 bounds for sums of multiple stochastic integrals. 
 In Section~\ref{s3} we show that graph counts 
 can be represented as sums of multiple stochastic
 integrals, and we derive Kolmogorov distance bounds
 for the renormalized count of subgraphs in $\mathbb{G}_n(p_n)$ that are isomorphic to
 a fixed graph.  
 \section{Notation and preliminaries}
 \label{s1}
 In this section we recall some background notation and results
 on the stochastic analysis of Bernoulli processes,
 see \cite{prsurvey} for details.
 Consider a sequence
 $(X_n)_{n\in \inte}$ of
 independent  identically distributed Bernoulli random variables
 with
 $P(X_n=1) =p$ and $P(X_n = -1)=q$,
 $n\in \inte$,
 built as the sequence of
 canonical projections on $\Omega := \{-1,1\}^{\mathbb{N}}$.
 For any $F:\Omega\to\R$ we consider the
 $L^2(\Omega\times\mathbb{N})$-valued 
 finite difference operator $D$ defined
 for any  $\omega=(\omega_0,\omega_1,\ldots)\in\Omega$
 by 
 \begin{equation}
   \label{fdb} 
D_k F(\omega)=\sqrt{pq}(F(\omega_{+}^k)-F(\omega_{-}^{k})),\quad k\in\mathbb{N}, 
\end{equation} 
 where 
 we let
$$
\omega_{+}^k:=(\omega_0,\ldots,\omega_{k-1},+1,\omega_{k+1},\ldots)
\mbox{~~and~~} 
\omega_{-}^k:=(\omega_0,\ldots,\omega_{k-1},-1,\omega_{k+1},\ldots),
\quad k \in \inte, 
$$
 and $DF := ( D_k F )_{k\in \inte}$.
 The $L^2$ domain of
$D$ is given by
$$
\mathrm{Dom}(D)
 =\{F\in
L^2(\Omega):\,\,\mathrm{E}[\|DF\|_{\ell^2(\mathbb{N})}^2]<\infty\}.
$$
 We let $( Y_n )_{n\geq 0}$ denote the sequence of
 centered and normalized random variables defined by
$$
Y_n :=\frac{q-p+X_n}{2\sqrt{p q}}, \qquad n \in \inte. 
$$
 Given $n\geq 1$, we denote by $\ell^2(\mathbb{N})^{\otimes
n}=\ell^2(\mathbb{N}^{n})$ the class of square-summable functions on
$\mathbb{N}^n$, we denote by $\ell^2(\mathbb{N})^{\circ n}$
the subspace of $\ell^2(\mathbb{N})^{\otimes n}$ formed by
functions that are symmetric in $n$ variables.
We let
$$
I_n(f_n )=\sum_{(i_1,\ldots,i_n)\in\Delta_n}f_n(i_1,\ldots,i_n)Y_{i_1}\cdots Y_{i_n}
$$
 denote the discrete multiple stochastic integral
 of order $n$ of $f_n$ in the subspace
 $\ell^2_{\mathfrak{s}} (\Delta_n )$ of $\ell^2(\mathbb{N})^{\circ n}$
 composed of symmetric kernels that vanish on
 diagonals, i.e. on the complement of
$$
\Delta_n=\{(k_1,\ldots,k_n)\in\mathbb{N}^n:\,\,k_i\neq
k_j,\,\,1\leq i<j\leq n\},\quad\text{$n\geq 1$}.
$$
The multiple stochastic integrals satisfy the isometry
and orthogonality relation 
 \begin{equation}
   \label{isomf} 
\mathrm{E}[I_n(f_n)I_m(g_m)]
=\ind_{\{n=m\}}
 n!\langle f_n,g_m\rangle_{\ell^2_{\mathfrak{s}} (\Delta_n )}
,
\end{equation} 
 $f_n \in\ell^2_{\mathfrak{s}} (\Delta_n )$, $g_m \in\ell^2_{\mathfrak{s}} (\Delta_m )$,
 cf. e.g. Proposition~1.3.2 of \cite{privaultbk2}.
The finite difference operator $D$ acts on multiple stochastic
integrals as follows:
\begin{equation}\nonumber 
D_k I_n(f_n )
=nI_{n-1}(f_n (*,k)\ind_{\Delta_n}(*,k))
=nI_{n-1}(f_n (*,k))
,
\end{equation}
$k\in\mathbb{N}$, $f_n \in\ell^2_{\mathfrak{s}} (\Delta_n )$,
and it satisfies the finite difference product rule
\begin{equation}
\label{prodrule} 
D_k(FG)  = FD_kG+GD_kF-\frac{X_k}{\sqrt{pq}} D_kFD_kG,
 \qquad k\in \inte. 
\end{equation} 
for $F,G:\Omega \rightarrow \real$,
see Propositions~7.3 and 7.8 of \cite{prsurvey}.
\\
 
Due to the chaos representation property of Bernoulli
random walks, any square integrable $F$
may be represented as $F=\sum_{n\geq 0}I_n(f_n)$,
$f_n\in\ell^2_{\mathfrak{s}} (\Delta_n )$, and the $L^2$ domain of
$D$ can be rewritten as
\begin{align}
\mathrm{Dom}(D)&=\left\{F=\sum_{n\geq 0}I_n(f_n) \ : \ \sum_{n\geq
1}n\,n!\|f_n\|_{\ell^2(\mathbb{N})^{\otimes
n}}^2<\infty\right\}.\nonumber
\end{align}
The Ornstein-Uhlenbeck operator $L$ is defined on the domain 
\begin{align}
\mathrm{Dom}(L) :=
 \left\{
 F=\sum_{n\geq 0}I_n(f_n) \ : \ \sum_{n\geq
1}n^2\,n!\|f_n\|_{\ell^2(\mathbb{N})^{\otimes
n}}^2<\infty\right\}
\nonumber
\end{align}
by 
$$
LF=-\sum_{n=1}^\infty nI_n(f_n).
$$
The inverse of $L$, denoted by $L^{-1}$, is defined on
the subspace of $L^2(\Omega)$ composed of centered random variables
by 
$$
L^{-1}F=-\sum_{n = 1}^\infty \frac{1}{n} I_n(f_n),
$$
 with the convention $L^{-1} F = L^{-1} ( F - \mathrm{E} [F] )$ in case
 $F$ is not centered. 
 Using this convention,
 the duality relation \eqref{prdual}
 shows that for any $F,G\in\mathrm{Dom}(D)$ 
 we have the covariance identity 
\begin{equation}\label{eq:covariance2}
\mathrm{Cov}(F,G)=\mathrm{E}[
 G ( F - \mathrm{E} [ F ] ) ]=\mathrm{E}\left[\langle
DG,-DL^{-1}F\rangle_{\ell^2(\mathbb N)}\right]. 
\end{equation}
 The divergence operator $\delta$ is 
 the linear mapping defined  as
$$
 \delta ( u ) 
 = 
 \delta (I_{n}(f_{n+1}( * ,\cdot )))=I_{n+1}(\tilde{f}_{n+1}), \quad 
 f_{n+1}\in \ell^2_{\mathfrak{s}} (\Delta_n ) \otimes \ell^{2}(\inte),
$$
 for $(u_k)_{k\in\inte}$ of the form 
$$ 
 u_k = 
 I_n (f_{n+1}( * , k )), 
 \qquad 
 k\in\inte 
, 
$$ 
 in the space 
$${\cal U} = \left\{
 \sum_{k=0}^n 
 I_k (f_{k+1} (*,\cdot ) ),
 \quad
 f_{k+1} \in \ell^2_{\mathfrak{s}} (\Delta_k ) \otimes \ell^2 (\inte), \
 \ k=, n \in\inte \right\} 
 \subset L^2 (\Omega \times \inte )
 $$ 
 of finite sums of multiple integral processes, 
 where $\tilde{f}_{n+1}$ denotes the symmetrization of $f_{n+1}$
 in $n+1$ variables, i.e. 
$$\tilde{f}_{n+1} (k_1,\ldots ,k_{n+1}) 
 = \frac{1}{n+1}
 \sum_{i=1}^{n+1}
 f_{n+1} (k_1,\ldots ,k_{k-1},k_{k+1}, \ldots ,k_{n+1},
 k_i ) 
. 
$$ 
 The operators $D$ and $\delta$ are closable
  with respective domains $\Dom ( D)$ and $\Dom ( \delta )$, 
  built as the completions of ${\cal S}$ and ${\cal U}$, 
  and they satisfy the duality relation  
 \begin{equation}
   \label{prdual}
   \E [\langle DF,u\rangle_{\ell^2 (\inte )} ]
 = \E [F\delta (u) ], \quad F\in \Dom ( D), \ u\in \Dom ( \delta ), 
 \end{equation}
 see e.g. Proposition~9.2 in \cite{prsurvey},
 and the isometry property
\begin{eqnarray} 
\nonumber 
 \E [ | \delta (u)| ^2 ]
 & = & \E [\Vert u \Vert_{\ell^2(\inte )}^2 ]
 + \E \Bigg[ \sum_{k,l=0\atop k \not= l}^\infty D_ku_l D_lu_k
  - \sum_{k=0}^\infty ( D_ku_k)^2 
  \Bigg]
 \\
 \label{skois}
 & \leq & \E [\Vert u \Vert_{\ell^2(\inte )}^2 ]
 + \E \Bigg[ \sum_{k,l=0\atop k \not= l}^\infty D_ku_l D_lu_k \Bigg], 
 \qquad u \in {\cal U}, 
\end{eqnarray} 
cf. Proposition~9.3 of \cite{prsurvey} and Satz~6.7 in \cite{mantei}.
 Letting $(P_t)_{t\in \real_+} = (e^{tL})_{t\in \real_+}$ 
 denote the Orsntein-Uhlenbeck semi-group defined as
$$P_t F = \sum_{n=0}^\infty e^{-nt} I_n(f_n),
 \qquad t\in \real_+,$$
 on random variables $F\in L^2(\Omega )$ of the form  
 $\displaystyle F = \sum_{n=0}^\infty I_n(f_n)$, 
 the Mehler formula states that 
\begin{equation} 
\label{istheou} 
 P_t F = \E [F(X (t) ) \mid X (0) ], 
 \qquad t\in \real_+, 
\end{equation}
where $(X (t) )_{t\in \real_+}$ is the Ornstein-Uhlenbeck
process associated to the semi-group
$(P_t)_{t\in \real_+}$,
cf. Proposition~10.8 of \cite{prsurvey}. 
As a consequence of the representation
\eqref{istheou} of $P_t$ we can deduce the bound
\begin{equation}
  \label{mehler} 
\E [ | D_k L^{-1} F|^\alpha ] \leq 
\E [ | D_k F|^\alpha ]
,
\end{equation}
for every $F\in \Dom (D)$ and $\alpha \geq 1$,
see Proposition~3.3 of \cite{reichenbachsAoP}.
The following Proposition~\ref{prop:d_Kdescrete} is
a consequence of
Proposition~4.1 in \cite{reichenbachsAoP},
see also Theorem~3.1 in \cite{reichenbachs}.
\begin{prop}\label{prop:d_Kdescrete}
For $F\in \mathrm{Dom}(D)$ with $\E[F]=0$ we have
\begin{align*}
  d_K(F,\mathcal{N} )\leq & |1-\E[F^2]|
  + \sqrt{ \Var[\langle DF,-DL^{-1}F\rangle_{\ell^2(\N)}] }
  \\
& +\frac{1}{2\sqrt{pq}} \sqrt{\sum_{k=0}^\infty\E[(D_kF)^4]}\(
\sqrt{\E\big[ F^2 \big]}
+\sqrt{\sum_{k=0}^\infty\E[(FD_kL^{-1}F)^2]}\)\\
&+ \frac{1}{\sqrt{pq}}
\sup_{x\in \R} \E[\langle D \mathbf1_{\{F>x\}},DF|DL^{-1}F|\rangle_{\ell^2(\N)}].
\end{align*}
\end{prop}
\begin{Proof}
 By Proposition~4.1 in \cite{reichenbachsAoP} we have 
 \begin{eqnarray} 
   \nonumber 
   d_K(F,\mathcal{N} ) & \leq & \E[|1-\langle DF,-DL^{-1}F\rangle_{\ell^2(\N)}|]
   \\
   \label{eq:dK1-1}
 & &    +
   \frac{\sqrt{2\pi}}{8} (pq)^{-1/2} \E[\langle | DF |^2,|DL^{-1}F|\rangle_{\ell^2(\N)}]
   \qquad 
   \\
   \label{eq:dK1-2}
   & & +\frac12(pq)^{-1/2}\E[\langle | DF |^2,|F DL^{-1}F|\rangle_{\ell^2(\N)}]
   \\
\nonumber 
& & + (pq)^{-1/2}\sup_{x\in \R} \E[\langle D \mathbf1_{\{F>x\}}, DF |DL^{-1}F|\rangle_{\ell^2(\N)}]. 
\end{eqnarray} 
On the other hand, the covariance identity \eqref{eq:covariance2}
 shows that 
$\E[ | \langle DF,-DL^{-1}F\rangle_{l^2(\N)} | ]=\Var F$,
hence by the Cauchy-Schwarz and triangular inequalities we get 
\begin{eqnarray*}
  \lefteqn{
    \! \! \! \! \! \! \! \! \! \! \! \! \!
    \mathrm{E}\left[\Big|1-\langle D F,-D L^{-1}F\rangle_{\ell^2(\mathbb{N})}\Big|\right]
 \leq 
\Big\|1-\langle D F,-D L^{-1}F\rangle_{\ell^2(\mathbb{N})}\Big\|_{L^2(\Omega)}
  }
  \\
& \leq &
 | 1
 -
 \Vert
 F
 \Vert_{L^2(\Omega)}^2
 |
 +
 \Vert
 \langle D F,-D L^{-1}F\rangle_{\ell^2(\mathbb{N})}
 -
 \Vert
 F
 \Vert_{L^2(\Omega)}^2
 \Vert_{L^2(\Omega)}
 \\
 & = &
 |1-\Var[F]|+
 \sqrt{ \Var[\langle DF,-DL^{-1}F\rangle_{\ell^2(\N)} ]}. 
\end{eqnarray*}
 Next, we have 
\begin{align*}
 \E \big[\Vert
 D L^{-1} I_n(f_n) \Vert_{\ell^2(\N)}^2 \big]
 &=\sum_{k=0}^\infty\E[(I_{n-1} (f_n(k,\cdot)))^2]
\\
 & = (n-1)!
 \sum_{k=0}^\infty \|f_n(k,\cdot)\|_{\ell^2(\N)^{\otimes (n-1)}}^2
\\
&=(n-1)!\|f_n\|_{\ell^2(\N)^{\otimes n}}^2
  \\
& \leq n!\|f_n\|_{\ell^2(\N)^{\otimes n}}^2
\\
& = \E\[|I_n(f_n)|^2\],
\end{align*}
and consequently, by the orthogonality relation
\eqref{isomf} we have
$$
\E \big[\Vert DL^{-1} F \Vert_{\ell^2(\N)}^2 \big]
\leq \E\big[F^2\big]
$$
for every $F\in L^2(\Omega)$, 
hence \eqref{eq:dK1-1} is bounded by 
\begin{align*}
  \E[ \langle |D L^{-1}F| , | DF |^2 \rangle_{\ell^2(\inte )} ]
  &\leq
  \E \left[
    \sqrt{
      \sum_{k=0}^\infty |D_kL^{-1}F|^2
      \sum_{k=0}^\infty |D_kF|^4
    }
  \right]\\
  &\leq
  \sqrt{     \E\left[ \sum_{k=0}^\infty |D_kL^{-1}F|^2 \right]}
  \sqrt{\E \left[ \sum_{k=0}^\infty  (D_kF)^4 \right]}
\\
  & = 
  \sqrt{ \E \big[\Vert DL^{-1} F \Vert_{\ell^2(\N)}^2 \big] }
  \sqrt{\E \left[ \sum_{k=0}^\infty  (D_kF)^4 \right]}
  \\
  &\leq
  \sqrt{\E[F^2]} \sqrt{\E \left[ \sum_{k=0}^\infty  (D_kF)^4 \right]}. 
\end{align*}
Eventually, regarding the third term \eqref{eq:dK1-2}, 
by the Cauchy-Schwarz inequality we find 
\begin{align*}
  \E\big[\langle
    (DF)^2,|F DL^{-1}F|\rangle_{\ell^2(\N)}\big]\leq
  \sqrt{
  \sum_{k=0}^\infty\E\big[(D_kF)^4\big]}\sqrt{\sum_{k=0}^\infty\E \big[(FD_kL^{-1}F)^2\big]}. 
\end{align*}
\end{Proof}
Finally, given $f_n\in\ell^2_{\mathfrak{s}} (\Delta_n )$ and
$g_m\in\ell^2_{\mathfrak{s}} (\Delta_m )$
we have the multiplication formula 
\begin{equation}
\label{e8}
I_n(f_n)I_m(g_m)=\sum_{s=0}^{2 \min (n , m)}I_{n+m-s}(h_{n,m,s}) ,
\end{equation}
see Proposition~5.1 of \cite{privaulttorrisi4},
provided that the functions
$$
 h_{n,m,s} : =
 \sum_{s\leq 2i \leq 2 \min (s , n , m)} i!
 \binom{n}{i}
 \binom{m}{i}
 \binom{i}{s-i}
 \left(
 \frac{q-p}{2\sqrt{pq}}
 \right)^{2i-s}
 f_n \tilde{\star}_i^{s-i} g_m
$$
belong to $\ell^2_{\mathfrak{s}} (\Delta_{n+m-s})$, $0\leq s \leq
2 \min (n , m)$,
 where $f_n \tilde{\star}_k^l g_m$ 
 is defined as the symmetrization in $n+m-k-l$ variables
 of the contraction $f_n \star_k^l g_m$ 
 defined as 
\begin{eqnarray*}
  \lefteqn{
    \! \! \! \! \! \! \! \! \! \! \! \! \! \! \! \! \! \!
    \! \! \! \! \! \! \! \! \! \! \! \! \! \! \! \! \! \!
    f_n \star_k^l g_m
    ( a_{l+1},\ldots ,a_n,b_{k+1},\ldots ,b_m) = 
\ind_{\Delta_{n+m-k-l}} ( a_{l+1},\ldots ,a_n,b_{k+1},\ldots ,b_m )
  }
  \\
  & & \times
  \sum_{a_1,\ldots ,a_l \in \mathbb{N}}f_n(a_1,\ldots ,a_n)g_m(a_1,\ldots ,a_k,b_{k+1} , \ldots , b_m), 
\end{eqnarray*}
 $0\leq l\leq k$,
and the symbol $\sum_{s\leq 2i \leq 2 \min (s , n , m)}$
means that the sum is taken over all the integers $i$
in the interval $[s/2, \min ( s , n , m )]$. 
We close this section with the following
Proposition~\ref{prop:fg<f^2+g^2}.
\begin{prop}\label{prop:fg<f^2+g^2}
 Let $f_n \in \ell^2_{\mathfrak{s}} (\Delta_n )$ and $g_m \in \ell^2_{\mathfrak{s}} (\Delta_m )$
 be symmetric functions. 
 For  $0\leq l< k\leq \min ( n,m )$ we have 
\begin{align}\label{eq:l<k} 
\left\| f_{n}\star_k^{l} g_m \right\|^2_{\ell^2(\N)^{\otimes (m+n-k-l)}}&\leq \frac12\left\| f_{n}\star_{n}^{l+n-k} f_{n}\right\|^2_{\ell^2(\N)^{\otimes (k-l)}}+\frac12\left\| g_m \star_{m}^{l+m-k} g_m \right\|^2_{\ell^2(\N)^{\otimes (k-l)}},
\end{align}
and 
\begin{align}\label{eq:l=k}
  \left\| f_{n}\star_k^{k} g_m \right\|^2_{\ell^2(\N)^{\otimes (m+n-2k)}}&\leq
  \frac{1}{2}
  \left\| f_{n}\star_{n-k}^{n-k} f_{n}\right\|^2_{\ell^2(\N)^{\otimes 2k}}
  +
  \frac{1}{2}
  \left\| g_m \star_{m-k}^{m-k} f_{m}\right\|^2_{\ell^2(\N)^{\otimes 2k}}.
\end{align}
\end{prop}
\begin{Proof} H\"older's inequality applied twice gives us
\begin{eqnarray*}
  \lefteqn{
    \! \! \! \! \! \! \! \! \! \!
    \left\| f_{n}\star_k^{l} g_m \right\|^2_{\ell^2(\N)^{\otimes (m+n-k-l)}}
    =\sum_{z_1\in\N^{n-k}}\sum_{z_2\in\N^{m-k}}\sum_{y\in\N^{k-l}}\(\sum_{x\in\N^l}f_n(x,y,z_1)g_m(x,y,z_2)\)^2
  }
  \\
  & \leq & \sum_{y\in\N^{k-l}}\sum_{z_1\in\N^{n-k}}\sum_{z_2\in\N^{m-k}}\(\sum_{x\in\N^l}f_n^2(x,y,z_1)\sum_{x\in\N^l}g_m^2(x,y,z_2)\)
  \\
  &\leq & 
  \sqrt{
    \sum_{y\in\N^{k-l}}\(\sum_{z_1\in\N^{n-k}}\sum_{x\in\N^l}f_n^2(x,y,z_1)\)^2\sum_{y\in\N^{k-l}}\(\sum_{z_1\in\N^{m-k}}\sum_{x\in\N^l}g_m^2(x,y,z_2)\)^2 }
  \\[6pt]
  &= & 
  \left\| f_{n}\star_{n}^{l+n-k} f_{n}\right\|_{\ell^2(\N)^{\otimes (k-l)}}
  \left\| g_m \star_{m}^{l+m-k} g_m \right\|_{\ell^2(\N)^{\otimes (k-l)}}
    \\[6pt]
  & \leq &  
    \frac{1}{2}
    \left\| f_{n}\star_{n}^{l+n-k} f_{n}\right\|^2_{\ell^2(\N)^{\otimes (k-l)}}
    +
    \frac{1}{2}
    \left\| g_m \star_{m}^{l+m-k} g_m \right\|^2_{\ell^2(\N)^{\otimes (k-l)}}. 
\end{eqnarray*}
To derive the second assertion, we proceed as follows: 
\begin{align*}
&\left\| f_{n}\star_k^{k} g_m \right\|^2_{\ell^2(\N)^{\otimes (m+n-2k)}}\\
&=\sum_{y\in\N^{n-k}}\sum_{z\in\N^{m-k}}\sum_{x_1\in\N^k}\sum_{x_2\in\N^k}f_n(x_1,y)g_m(x_1,z)f_n(x_2,y)g_m(x_2,z)\\
&=\sum_{x_1\in\N^k}\sum_{x_2\in\R^k}\(\sum_{y\in\N^{n-k}}f_n(x_1,y)f_n(x_2,y)\)\(\sum_{z\in\N^{m-k}}g_m(x_1,z)g_m(x_2,z)\)\\
  &\leq
  \frac{1}{2}
  \sum_{x_1\in\N^k}\sum_{x_2\in\R^k}\(\sum_{ y\in \N^{n-k}}f_n(x_1,y)f_n(x_2,y)\)^2
  +
  \frac{1}{2}
  \sum_{x_1\in\N^k}\sum_{x_2\in\R^k}\(\sum_{ z\in \N^{m-k}}g_m(x_1,z)g_m(x_2,z)\)^2
  \\
 & = 
    \frac{1}{2}
  \left\| f_{n}\star_{n-k}^{n-k} f_{n}\right\|^2_{\ell^2(\N)^{\otimes 2k}}
  +
  \frac{1}{2}
  \left\| g_m \star_{m-k}^{m-k} f_{m}\right\|^2_{\ell^2(\N)^{\otimes 2k}}. 
\end{align*}
\end{Proof} 
\section{Kolmogorov bounds for sums of multiple stochastic integrals} 
\label{s2}
Wasserstein bounds have been obtained for discrete
multiple stochastic integrals in Theorem~4.1 of \cite{nourdin3}
in the symmetric case $p=q$
and in Theorems~5.3-5.5 of \cite{privaulttorrisi4} 
in the possibly nonsymmetric case,
and have been extended to the Kolmogorov distance
in the symmetric case $p=q$ in Theorem~4.2 of \cite{reichenbachs}.
The following result provides a Kolmogorov distance bound which
further extends Theorem~4.2 of \cite{reichenbachs}
from multiple stochastic integrals to
sums of multiple stochastic integrals in the nonsymmetric case.
\begin{theorem}\label{prop:dKsumI}
 For any finite sum 
 $$ F=\sum_{k=1}^nI_k(f_k)$$
 of discrete multiple stochastic integrals
 with $f_k \in \ell^2_{\mathfrak{s}} (\Delta_k )$, $k=1,\ldots, n$, we have
\begin{align*}
d_K(F,\mathcal{N} )\leq & C_n \big( |1-\Var[F]|+\sqrt{R_F} \big),
\end{align*}
for some constant $C_n>0$ depending only on $n$, where
\begin{equation} 
\label{sjks}
 R_F:= \sum_{0\leq l< i\leq n}(pq)^{l-i}\left\| f_{i}\star_{i}^{l} f_{i}\right\|^2_{\ell^2(\N)^{\otimes (i-l)}}+\sum_{1\leq l< i\leq n}
 \(
 \left\| f_{l}\star_{l}^{l} f_{i}\right\|^2_{\ell^2(\N)^{\otimes (i-l)}}
 +
 \left\| f_{i}\star_{l}^{l} f_{i}\right\|^2_{\ell^2(\N)^{\otimes 2(i-l)}}
 \). 
\end{equation} 
\end{theorem}
\begin{Proof}
 We introduce
$$R'_F:=\sum_{1\leq i\leq j\leq n}\sum_{k=1}^{i}\sum_{l=0}^k\mathbf1_{\{i=j=k=l\}^c}
(pq)^{l-k}
\left\| f_{i}\star_k^{l} f_{j}\right\|^2_{\ell^2(\N)^{\otimes (i+j-k-l)}}.
$$
Since it holds that $R'_F\lesssim R_F$, it is enough
to prove the required inequality with $R'_F$ instead of $R_F$.
Indeed, by the inequality \eqref{eq:l<k},
all the components of $R'_F$ for $0\leq l<k\leq i,j$, are dominated by those for $0\leq l<k=i=j$, and also, by the inequality \eqref{eq:l=k}, 
the ones where $1\leq k=l< i \leq j$,
are dominated by the components where $1\leq l=k<i=j$. Finally, the components for $1\leq k=l=i<j$ remain unchanged. 
\\
 
We will estimate components in the inequality from Proposition~\ref{prop:d_Kdescrete}.
We have
\begin{align*}
  D_r F=(i+1)\sum_{i=0}^{n-1}I_i \(f_{i+1}(r,\cdot)\),\quad
  \mbox{and} \quad 
D_r L^{-1}F = \sum_{i=0}^{n-1}I_i \(f_{i+1}(r,\cdot)\),
\qquad r \in\inte, 
\end{align*}
hence by the multiplication formula \eqref{e8} we find 
\begin{equation}
  \label{XX}
  (D_r F)^2=\sum_{0\leq i\leq j\leq n-1}\sum_{k=0}^{i}\sum_{l=0}^kc_{i,j,l,k}\(\frac{q-p}{\sqrt{pq}}\)^{k-l}I_{i+j-k-l}\(f_{i+1}(r,\cdot)\tilde{\star}_k^lf_{j+1}(r,\cdot)\)
\end{equation} 
and
\begin{equation}
  \label{XLX}
D_r FD_r L^{-1}F=\sum_{0\leq i\leq j\leq n-1}\sum_{k=0}^{i}\sum_{l=0}^kd_{i,j,l,k}\(\frac{q-p}{\sqrt{pq}}\)^{k-l}I_{i+j-k-l}\(f_{i+1}(r,\cdot)\tilde{\star}_k^lf_{j+1}(r,\cdot)\),
\end{equation}
for some $c_{i,j,l,k}$, $d_{i,j,l,k} \geq0$.
Applying the isometry relation \eqref{isomf} to \eqref{XX} and
using the bound $\| \tilde{f}_n\|_{\ell^2(\N )^{\otimes n}} \leq \| f_n \|_{\ell^2(\N )^{\otimes n}}$, 
$f_n \in \ell^2(\N )^{\otimes n}$, we get, 
writing $f\lesssim g$ whenever $f< C_n g$ for some 
 constant $C_n >0$ depending only on $n$, 
\begin{align}\nonumber
  \sum_{r=0}^\infty  \E \left[ | D_r F |^4  \right]&\lesssim \sum_{0\leq i\leq j\leq n-1}\sum_{k=0}^{i}\sum_{l=0}^k\sum_{r=0}^\infty \(\frac{q-p}{\sqrt{pq}}\)^{2k-2l}\left\|f_{i+1}(r,\cdot)\star_k^lf_{j+1}(r,\cdot)\right\|^2_{\ell^2(\N)^{\otimes (i+j-k-l)}}
  \\
  \nonumber
&= \sum_{0\leq i\leq j\leq n-1}\sum_{k=0}^{i}\sum_{l=0}^k \(\frac{q-p}{\sqrt{pq}}\)^{2k-2l}\left\|f_{i+1} \star_{k+1}^lf_{j+1} \right\|^2_{\ell^2(\N)^{\otimes (i+j-k-l+1)}}
\\
\nonumber
&= \sum_{1\leq i\leq j\leq n}\sum_{k=1}^{i}\sum_{l=0}^{k-1} \(\frac{q-p}{\sqrt{pq}}\)^{2k-2l-2}\left\|f_{i} \star_k^lf_{j} \right\|^2_{\ell^2(\N)^{\otimes (i+j-k-l)}}
\\
\label{aux1}
&\leq pqR'_F. 
\end{align}
Furthermore, by \eqref{XLX} it follows that 
\begin{align*}
  &\langle D F, D L^{-1}F\rangle-\E\[\langle D F, D L^{-1}F\rangle\]
  \\
&=\sum_{r=0}^\infty\sum_{0\leq i\leq j\leq n-1}\sum_{k=0}^{i}\sum_{l=0}^kc_{i,j,l,k}\mathbf1_{\{i=j=k=l\}^c}\(\frac{q-p}{\sqrt{pq}}\)^{k-l}I_{i+j-k-l}\(f_{i+1}(r,\cdot)\tilde{\star}_k^lf_{j+1}(r,\cdot)\)\\
&=\sum_{0\leq i\leq j\leq n-1}\sum_{k=0}^{i}\sum_{l=0}^kc_{i,j,l,k}\mathbf1_{\{i=j=k=l\}^c}\(\frac{q-p}{\sqrt{pq}}\)^{k-l}I_{i+j-k-l}\(\sum_{r=0}^\infty f_{i+1}(r,\cdot)\tilde{\star}_k^lf_{j+1}(r,\cdot)\)\\
&=\sum_{0\leq i\leq j\leq n-1}\sum_{k=0}^{i}\sum_{l=0}^kc_{i,j,l,k}\mathbf1_{\{i=j=k=l\}^c}\(\frac{q-p}{\sqrt{pq}}\)^{k-l}I_{i+j-k-l}\( f_{i+1} \tilde{\star}_{k+1}^{l+1}f_{j+1} \), 
\end{align*}
thus we get
\begin{eqnarray}
  \nonumber
  \lefteqn{
    \! \! \! \! \! \! \! \! \! \! \! \! \! \! \! \! \! \! \! \! \! \! \! \! \! \! \! \! \! \! \! \! \! \! \! \! \! \! \! 
    \Var\[\langle D F, -D L^{-1}F\rangle\]
    \lesssim\sum_{0\leq i\leq j\leq n-1}\sum_{k=0}^{i}\sum_{l=0}^k\mathbf1_{\{i=j=k=l\}^c}\(\frac{1}{\sqrt{pq}}\)^{2k-2l}\left\| f_{i+1}\star_{k+1}^{l+1} f_{j+1}\right\|^2_{\ell^2(\N)^{\otimes (i+j-k-l)}}
  }
  \\\nonumber
&= & \sum_{1\leq i\leq j\leq n}\sum_{k=1}^{i}\sum_{l=1}^k\mathbf1_{\{i=j=k=l\}^c}\frac{1}{(pq)^{k-l}}\left\| f_{i}\star_k^{l} f_{j}\right\|^2_{\ell^2(\N)^{\otimes (i+j-k-l)}}\\
\label{eq:term1} 
&\leq & R'_F.
\end{eqnarray}
Next, we have 
\begin{align*}
  &\sum_{k=0}^\infty
  \E \big[ (FD_kL^{-1}F)^2 \big]=
  \E \left[
    F^2\sum_{k=0}^\infty(D_kL^{-1}F)^2
    \right]\leq \sqrt{\E\[F^4\]}
  \sqrt{ \E\[\(\sum_{k=0}^\infty(D_kL^{-1}F)^2\)^2\]} 
\end{align*}
 and \eqref{e8} and \eqref{isomf} show that 
\begin{align*}
  \E\[F^4\]&\lesssim \E\[\(\sum_{1\leq i\leq j\leq n}\sum_{k=0}^{i}\sum_{l=0}^k
  \left|
  \frac{q-p}{\sqrt{pq}}
  \right|^{k-l}I_{i+j-k-l}\( f_{i} \tilde{\star}_k^{l}f_{j} \)\)^2\] \\
&\lesssim \sum_{1\leq i\leq j\leq n}\sum_{k=0}^{i}\sum_{l=0}^k
(pq)^{l-k}
\| f_{i} \star_k^{l}f_{j} \|_{\ell^2(\N)^{\otimes  (i+j-k-l)}}^2
\\
&\lesssim R'_F
+
\sum_{i=1}^n \left\| f_{i}\star_{i}^{i} f_{i}\right\|^2_{\ell^2(\N)^{\otimes 0}}
+
\sum_{1\leq i < j\leq n} \| f_{i} \star_0^0 f_{j} \|_{\ell^2(\N)^{\otimes (i+j)}}^2
\\
& = R'_F
+
\sum_{i=1}^n \left\| f_i \right\|^4_{\ell^2(\N)^{\otimes i}}
+
\sum_{1\leq i < j\leq n} \| f_i \|_{\ell^2(\N)^{\otimes i}}^2
\| f_j \|_{\ell^2(\N)^{\otimes j}}^2
\\
&\lesssim R'_F+(\Var[F])^2,
\end{align*}
 while as in \eqref{XX} and \eqref{XLX} we have 
\begin{align*}
  &\E\[\(\sum_{k=0}^\infty(D_kL^{-1}F)^2\)^2\]
 \\
& =
  \E\[
  \(
  \sum_{k=0}^\infty
  \sum_{0\leq i\leq j\leq n-1}\sum_{k=0}^{i}\sum_{l=0}^k\tilde{d}_{i,j,l,k}\(\frac{q-p}{\sqrt{pq}}\)^{k-l}I_{i+j-k-l}\(f_{i+1}(k,\cdot)\tilde{\star}_k^lf_{j+1}(k,\cdot)\)
  \)^2 
\]
\\
 &\lesssim \sum_{0\leq i\leq j\leq n-1}\sum_{k=0}^{i}\sum_{l=0}^k
  (pq)^{l-k}
\| f_{i+1} \star_{k+1}^{l+1}f_{j+1} \|^2_{\ell^2(\N)^{\otimes (i+j-k-l)}}
\\
 &= \sum_{1\leq i\leq j\leq n}\sum_{k=1}^i\sum_{l=1}^k
  (pq)^{l-k}
\| f_i \star_k^lf_j \|^2_{\ell^2(\N)^{\otimes (i+j-k-l)}}
\\
  &\lesssim R'_F
+
\sum_{i=1}^n \left\| f_{i}\star_{i}^{i} f_{i}\right\|^2_{\ell^2(\N)^{\otimes 0}}
+
\sum_{1\leq i < j \leq n} \| f_{i} \star_0^0 f_{j} \|_{\ell^2(\N)^{\otimes (i+j)}}^2
\\
& = R'_F
+
\sum_{i=1}^n \left\| f_i \right\|^4_{\ell^2(\N)^{\otimes i}}
+
\sum_{1\leq i < j\leq n} \| f_i \|_{\ell^2(\N)^{\otimes i}}^2
\| f_j \|_{\ell^2(\N)^{\otimes j}}^2
\\
 & \lesssim R'_F +(\Var[F])^2,
\end{align*}
 hence we get
\begin{align}\label{eq:term2}\sum_{k=0}^\infty\E[(FD_k L^{-1} F)^2]\lesssim R'_F+(\Var[F])^2.
\end{align}
We now deal with the last component in Proposition~\ref{prop:d_Kdescrete} similarly as it is done in proof of Theorem~4.2 in \cite{reichenbachs}. Precisely, by the
integration by parts formula \eqref{prdual}
and the Cauchy-Schwarz inequality we have
\begin{align}
  \nonumber 
  \sup_{x\in \R} \E\[\langle D \mathbf1_{\{F>x\}},DF|DL^{-1}F|\rangle_{\ell^2(\N)}\]&=\sup_{x\in \R} \E\[\mathbf1_{\{F>x\}}\delta\(DF|DL^{-1}F|\)\]
\\
\label{dklsfd}
  &\leq
\sqrt{ \E\[\(\delta\(DF|DL^{-1}F|\)\)^2\] }. 
\end{align}
Then, by the bound \eqref{skois},
the Cauchy-Schwarz inequality and the consequence \eqref{mehler}
of Mehler's formula \eqref{istheou}, we have 
\begin{align*}
  &\E\big[ \(\delta\(DF|DL^{-1}F|\)\)^2 \big]
  \\
  &\leq \E \big[\|DF|DL^{-1}F|\|^2_{\ell^2(\N)}\big]+\E\[\sum_{k,l=0}^\infty
  \big| D_k\(D_lF|D_lL^{-1}F|\)D_l\(D_kF|D_kL^{-1}F|\) \big|
  \]
  \\
  &\leq
  \sqrt{
    \E\big[\|DF\|^4_{\ell^4(\N)}\big]\E\big[\|DL^{-1}F\|^4_{\ell^4(\N)}\big]
    } 
  +\E\[\sum_{k,l=0}^\infty \(D_k\(D_lF|D_lL^{-1}F|\)\)^2 \]
  \\
  &\leq\E \big[\|DF\|^4_{\ell^4(\N)}\big]+\sum_{k,l=0}^\infty\E\big[
    \(D_k\(D_lF|D_lL^{-1}F|\)\)^2 \big].
\end{align*}
The first term in the last expression in bounded by $pqR'_F$ as shown in
\eqref{aux1},
and it remains to estimate the last expectation. By the
product rule \eqref{prodrule} and the bound $|D_k| F ||\leq |D_k F |$
obtained from the definition \eqref{fdb}
of $D$ and the triangle inequality, we get
\begin{align}
  \nonumber
  &\E\big[\(D_r\(D_sF|D_sL^{-1}F|\)\)^2 \big]
  \\
\nonumber
&=\E\[\(\(D_rD_sF|D_sL^{-1}F|\)+\(D_sFD_r|D_sL^{-1}F|\)-\frac{X_r}{\sqrt{pq}}\(D_rD_sFD_r|D_sL^{-1}F|\)\)^2 \]
\\
\label{dksds}
&\lesssim \E\[\(D_rD_sF\)^2\(D_sL^{-1}F\)^2+\(D_sF\)\(D_rD_sL^{-1}F\)^2+\frac{1}{pq}\(D_rD_sF\)^2\(D_rD_sL^{-1}F\)^2 \], 
\end{align}
$r,s \in \inte$.
By the Cauchy-Schwarz inequality  we get
\begin{align*}
\sum_{r,s=0}^\infty\E\big[\(D_rD_sF\)^2\(D_sL^{-1}F\)^2\big]
& =\E\[\sum_{s=0}^\infty\(D_sL^{-1}F\)^2\sum_{r=0}^\infty\(D_rD_sF\)^2\]\\
&\leq
\sqrt{
  \E\[\sum_{s=0}^\infty\(D_sL^{-1}F\)^4\]\E\[\sum_{s=0}^\infty\(\sum_{r=0}^\infty\(D_rD_sF\)^2\)^2\]}.
\end{align*}
The term $\E\big[\sum_{s=0}^\infty\(D_sL^{-1}F\)^4\big]$ can
be bounded by $pq R'_F$ as in \eqref{aux1}. To estimate the other
term we use the multiplication formula \eqref{e8} as in \eqref{XX}
to obtain
\begin{align*}
&\E\[\sum_{s=0}^\infty\(\sum_{r=0}^\infty\(D_rD_sF\)^2\)^2\]\\
  &\lesssim \sum_{s=0}^\infty\E\[\(\sum_{r=0}^\infty\sum_{0\leq i\leq j\leq n-2}\sum_{k=0}^{i}\sum_{l=0}^k
  \left|
  \frac{q-p}{\sqrt{pq}}
  \right|^{k-l}I_{i+j-k-l}\(f_{i+2}(s,r,\cdot)\tilde{\star}_k^lf_{j+2}(s,r,\cdot)\)\)^2\]\\
  &=c \sum_{s=0}^\infty\E\[\(\sum_{0\leq i\leq j\leq n-2}\sum_{k=0}^{i}\sum_{l=0}^k
  \left|
  \frac{q-p}{\sqrt{pq}}
  \right|^{k-l}I_{i+j-k-l}\(f_{i+2}(s,\cdot)\tilde{\star}_{k+1}^{l+1}f_{j+2}(s,\cdot)\)\)^2\]
  \\
  &\lesssim  \sum_{s=0}^\infty\sum_{0\leq i\leq j\leq n-2}\sum_{k=0}^{i}\sum_{l=0}^k
  (pq)^{l-k}\|f_{i+2}(s,\cdot) \star_{k+1}^{l+1}f_{j+2}(s,\cdot) \|^2_{\ell^2(\N)^{\otimes (i+j-k-l)}}\\
  &=  \sum_{0\leq i\leq j\leq n-2}\sum_{k=0}^{i}\sum_{l=0}^k
  (pq)^{l-k}\|f_{i+2} \star_{k+2}^{l+1}f_{j+2} \|^2_{\ell^2(\N)^{\otimes (i+j-k-l+1)}}\\
  &=  \sum_{2\leq i\leq j\leq n}\sum_{k=2}^{i}\sum_{l=1}^{k-1}
  (pq)^{l+1-k}\|f_{i} \star_k^{l}f_{j} \|^2_{\ell^2 (\N)^{\otimes (i+j-k-l)}}
  \\
   & \leq pqR'_F.
\end{align*}
The term $\sum_{r,s=0}^\infty\E\big[(D_sF)^2(D_rD_sL^{-1}F)^2\big]$
from \eqref{dksds} is similarly bounded by $pqR'_F$.
Regarding the last term, we have 
$$ 
\sum_{r,s=0}^\infty\E\big[(D_rD_sF)^2(D_rD_sL^{-1}F)^2\big]
\leq 
  \sqrt{
  \sum_{r,s=0}^\infty\E\[\(D_rD_sF\)^4\]\sum_{r,s=0}^\infty
  \E\big[(D_rD_sL^{-1}F)^4\big] }. 
$$ 
  Using the multiplication formula \eqref{e8}, both sums inside the
  above square root can be estimated as 
\begin{align*}
  &\sum_{r,s=0}^\infty\E\[\(\sum_{0\leq i\leq j\leq n-2}\sum_{k=0}^{i}\sum_{l=0}^k
  \left|
  \frac{q-p}{\sqrt{pq}}
  \right|^{k-l}I_{i+j-k-l}\(f_{i+2}(s,r,\cdot)\tilde{\star}_k^lf_{j+2}(s,r,\cdot)\)\)^2\]\\
 &\lesssim\sum_{r,s=0}^\infty\sum_{0\leq i\leq j\leq n-2}\sum_{k=0}^{i}\sum_{l=0}^k
 (pq)^{l-k}\|f_{i+2}(s,r,\cdot)\star_k^lf_{j+2}(s,r,\cdot)\|^2_{\ell^2(\N)^{\otimes (i+j-k-l)}}\\
 &=\sum_{0\leq i\leq j\leq n-2}\sum_{k=0}^{i}\sum_{l=0}^k
 (pq)^{l-k}\|f_{i+2} \star_{k+2}^{l}f_{j+2} \|^2_{\ell^2(\N)^{\otimes (i+j-k-l + 2)}}
\\
 &=\sum_{2\leq i\leq j\leq n}\sum_{k=2}^{i}\sum_{l=0}^{k-2}
 (pq)^{l+2-k}\|f_{i} \star_k^{l}f_{j} \|^2_{\ell^2(\N)^{\otimes (i+j-k-l)}}
\\
 & \lesssim (pq)^2R'_F. 
\end{align*}
 Combining this together we get 
$$
 \sum_{r,s=0}^\infty
 \E\big[ \big(D_r\big(D_sF|D_sL^{-1}F|\big)\big)^2 \big]\lesssim pq R'_F.
$$
and consequently, by \eqref{dklsfd} we find 
\begin{align}\label{eq:term3}
\sup_{x\in \R} \E\[\langle D \mathbf1_{\{F>x\}},DF|DL^{-1}F|\rangle_{\ell^2(\N)}\]\lesssim pq R'_F.
\end{align}
Applying \eqref{aux1}-\eqref{eq:term2} and \eqref{eq:term3}
to Proposition~\ref{prop:d_Kdescrete}, we get 
$$ 
  d_K(F,\mathcal{N} ) \lesssim |1-\Var[F]|+\sqrt{R'_F}
  \big(1+\Var[F]
  +\sqrt{\Var[F]}
  +\sqrt{R'_F}\big).
$$ 
If $R'_F\geq1$, or if $R'_F\leq1$ and $\Var[F]\geq2$, it is clear that  $d_K(F,\mathcal{N} )\lesssim |1-\Var[F]|+\sqrt{R'_F}$ since $d_K(F,\mathcal{N} )\leq 1$ by definition.
If $R'_F\leq1$ and  $\Var[F]\leq2$, we estimate $\Var[F]+\sqrt{\Var[F]}+\sqrt{R'_F}$ by a constant and also get the required bound. 
\end{Proof}
\section{Application to random graphs}
\label{s3}
In the sequel fix a
numbering $(1,\ldots , e_G)$ of the
edges in $G$ and 
we denote by $E_G \subset \inte^{e_G}$
the set of sequences of (distinct) edges that create a graph isomorphic to $G$,
i.e. a sequence $( e_{k_1} ,\ldots , e_{k_{e_G}} )$ belongs to $E_G$
if and only if the graph created by edges $e_{k_1},\ldots ,e_{k_{e_G}}$
is isomorphic to $G$.
The next lemma allows us to represent the number
of subgraphs as a sum of multiple stochastic integrals,
using the notation $P(X_k=1)=p$, $P(X_k=-1)=1-p=q$, $k\in\inte$. 
\begin{lemma}
  \label{lemma4.1}
  We have the identity 
\begin{align}\label{eq:tildeNasI}
\tilde{N}_G=\frac{N_G-\E[N_G]}{\sqrt{\Var[N_G]}}=\sum_{k=1}^{e_G}I_k(f_k),
\end{align}
where
$$
f_k(b_1,\ldots , b_k) := \frac{q^{k/2}p^{e_G - k/2}}{(e_G-k)!k!\sqrt{\Var[N_G]}}\,
\sum_{(a_1,\ldots , a_{e_{\scaleto{G}{3pt}} -k} ) \in
  \N^{e_{\scaleto{G}{3pt}} -k}} \mathbf1_{(a_1,\ldots , a_{e_{\scaleto{G}{3pt}} -k} , b_1,\ldots , b_k )\in E_G}.
$$
\end{lemma}
\begin{Proof} 
  We have
  \begin{align}\nonumber
  &  N_G =\frac1{e_G!2^{e_G}}\sum_{ b_1,\ldots ,b_{e_G} \in \inte}
    \mathbf1_{(b_1,\ldots ,b_{e_G})\in E_G}(X_{b_1}+1)\cdots (X_{b_{e_G}}+1)\\\nonumber
&=\frac1{e_G!2^{e_G}}\sum_{m=0}^{e_G}\({{e_G}\atop m}\)\sum_{b_1,\ldots ,b_{m} \in \inte}
    g_m(b_1,\ldots ,b_{m})X_{b_1}\cdots X_{b_{m}}
    \\\nonumber
    &=\frac1{e_G!2^{e_G}}\sum_{m=0}^{e_G}\({{e_G}\atop m}\)\sum_{k=0}^m\({m\atop k}\)
(p-q)^{m-k}
    \sum_{b_1,\ldots,b_k \in \inte } g_k(b_1,\ldots,b_k)
( X_{b_1} + q-p ) 
\cdots 
( X_{b_k} + q-p ) 
    \\\nonumber
    &=\frac1{e_G!2^{e_G}}\sum_{m=0}^{e_G}\({{e_G}\atop m}\)\sum_{k=0}^m\({m\atop k}\)I_k(g_k)(2\sqrt{pq})^k(p-q)^{m-k}
    \\
    \nonumber
    &=\frac1{e_G!2^{e_G}}\sum_{k=0}^{e_G}
\(e_G \atop k \)
(2\sqrt{pq})^kI_k(g_k)
    \sum_{m=k}^{e_G}\({{e_G-k}\atop {m-k}}\)(p-q)^{m-k}
    \\
    \nonumber
    &=\frac{1}{2^{e_G}}\sum_{k=0}^{e_G}
\frac{(2\sqrt{pq})^k}{(e_G-k)!k!} 
I_k(g_k)
(1+p-q)^{e_G-k}
\\
\nonumber
&=\sum_{k=0}^{e_G}
\frac{q^{k/2}p^{e_G-k/2}}{(e_G-k)!k!} 
I_k(g_k),
\end{align}
 where $g_k$ is the function defined as 
\begin{equation}
\label{gk}
g_k(b_1,\ldots , b_k)
:=\sum_{(a_1,\ldots , a_{e_G-k} ) \in\N^{e_G-k}}\mathbf1_{E_G}
(a_1,\ldots , a_{e_G-k} , b_1,\ldots , b_k ),
\quad (b_1,\ldots , b_k ) \in\N^k, 
\end{equation}
which shows \eqref{eq:tildeNasI} with
$$
f_k(b_1,\ldots , b_k ):=
\frac{q^{k/2}p^{e_G- k/2 }}{(e_G-k)!k! \sqrt{\Var[N_G]}}\,g_k(b_1,\ldots , b_k ).
$$ 
\end{Proof}
Next is the second main result of this paper. 
\begin{theorem}\label{thm:main}
Let  $G$ be a graph without isolated vertices. Then we have 
\begin{align*}
d_K (\tilde{N}_G,\mathcal{N} )&\lesssim\((1-p)\min_{\substack{ H\subset G\\e_H\geq1}}\big\{n^{v_H}p^{e_H}\big\}\)^{-1/2}.
\end{align*}
\end{theorem}
\begin{Proof}
 By \eqref{eq:tildeNasI} and Theorem~\ref{prop:dKsumI} we have
\begin{align}\label{eq:dKfromprop} 
d_K (\tilde{N}_G,\mathcal{N} )\lesssim  {\frac {\sqrt{R_G}}{\Var[N_G]}},
\end{align}
where, taking $g_k$ as in \eqref{gk}, by \eqref{sjks} we have  
\begin{align*}
  R_G = \ & \sum_{0\leq l< k\leq e_G}p^{4e_G-3k+l}q^{l+k}\left\| g_k\star_k^{l} g_k\right\|^2_{\ell^2(\N)^{\otimes (k-l)}}
 +\sum_{1\leq l< k\leq e_G}p^{4e_G-2k}q^{2k} \left\| g_k\star_{l}^{l} g_k\right\|^2_{\ell^2(\N)^{\otimes  2(k-l)}}
\\
  \ &
  +\sum_{1\leq l<k\leq e_G}p^{4e_G-l-k}q^{k+l}\left\| g_l \star_{l}^{l} g_k\right\|^2_{\ell^2(\N)^{\otimes (k-l)}}
\\
\leq \ & q \Big( \sum_{0\leq l< k\leq e_G}p^{4e_G-3k+l}\left\| g_k\star_k^{l} g_k\right\|^2_{\ell^2(\N)^{\otimes (k-l)}}
+\sum_{1\leq l<k\leq e_G}p^{4e_G-l-k}\left\| g_l \star_{l}^{l} g_k\right\|^2_{\ell^2(\N)^{\otimes (k-l)}}
\\
 \ &
+\sum_{1\leq l<k\leq e_G}p^{4e_G-2k}
\left\| g_k\star_{l}^{l} g_k\right\|^2_{\ell^2(\N)^{\otimes 2(k-l)}}
\Big)\\
= \ & (1-p) \big( S_1+S_2+S_3\big).
\end{align*}
It is now sufficient to show that 
\begin{align}\label{eq:SSS<}S_1+S_2+S_3\lesssim \max_{ H\subset G
    \atop
    e_H\geq1}n^{4v_G-3v_H}p^{4e_G-3e_H}.
\end{align}
Indeed, applying \eqref{varng} and \eqref{eq:SSS<} to \eqref{eq:dKfromprop} we get
\begin{align*}
  \frac{\sqrt{R_G}}{\Var[N_G]}&\lesssim \frac{\sqrt{1-p}
    \sqrt{\displaystyle
      \max_{ H\subset G \atop e_H\geq1}n^{4v_G-3v_H}p^{4e_G-3e_H}}}{(1-p)
    \displaystyle
    \max_{ H\subset G\atop e_H\geq1}n^{2v_G-v_H}p^{2e_G-e_H}}\\
  &=\frac{\(
    \displaystyle
    \min_{ H\subset G\atop e_H\geq1}n^{v_H}p^{e_H}\)^{-3/2}}{
    \sqrt{1-p} \(
    \displaystyle
    \min_{ H\subset G\atop e_H\geq1}n^{v_H}p^{e_H}\)^{-1}}\\
&=\((1-p)\min_{ H\subset G\atop e_H\geq1}n^{v_H}p^{e_H}\)^{-1/2}.
\end{align*}
 Thus
$$
 d_K (\tilde{N}_G,\mathcal{N} )\lesssim \frac{\sqrt{R_G}}{\Var[N_G]}\lesssim \((1-p)\min_{ H\subset G\atop e_H\geq1}n^{v_H}p^{e_H}\)^{-1/2}.
$$
 In order to estimate $S_1$, let us observe that 
\begin{align*}
  \left\| g_k \star_k^{l} g_k\right\|^2_{\ell^2(\N)^{\otimes (k-l)}}&=
  \sum_{a''\in\N^{k-l}}\(
  \sum_{a'\in\N^l}\(\sum_{a\in\N^{e_G-k}}\mathbf1_{E_G}\(a,a',a''\)\)^2\)^2
  \\
  &\approx\sum_{A \subset K_n
    \atop e_K=k-l}\(
  \sum_{ A \subset B \subset K_n \atop e_B =k}\(\sum_{
    B \subset G' \subset K_n \atop G'\sim G}1\)^2\)^2
  \\
  &\approx\sum_{K \subset G \atop e_K =k-l}n^{v_K}\(
  \sum_{K\subset H \subset G \atop e_H =k}n^{v_H-v_K}\(n^{v_G-v_H}\)^2\)^2
  \\[10pt]
  &\approx \max_{K\subset H \subset G
    \atop e_K=k-l, \ \! e_H =k}n^{4v_G-2v_H -v_K}.
\end{align*}
 Hence we have 
\begin{align*}
  S_1&\lesssim\sum_{0\leq l< k\leq e_G}p^{4e_G-3k+l}
  \max_{K\subset H \subset G
    \atop e_K=k-l, \ \! e_H =k}n^{4v_G-2v_H -v_K}\\
  &=\sum_{0\leq l< k\leq e_G}\max_{K\subset H \subset G
    \atop
    e_K=k-l, \ \! e_H =k}n^{4v_G-2v_H -v_K}p^{4e_G-2e_H -e_K}\\
&\lesssim \max_{K\subset H\subset G\atop e_K\geq1}n^{4v_G-2v_H -v_K}p^{4e_G-2e_H -e_K}.
\end{align*}
For a fixed $p$, let $H_0\subset G$, $e_{H_0}\geq1$, be the
subgraph of $G$ such that
\begin{equation}
\label{djlcn} 
n^{v_{H_0}}p^{e_{H_0}} = \min_{H \subset G, e_H\geq1}n^{v_H}p^{e_H}. 
\end{equation} 
 Then it is clear that
\begin{align}\nonumber
  S_1 & \lesssim
  \max_{K \subset H \subset G \atop
    e_K\geq1}n^{4v_G-2v_H -v_K}p^{4e_G-2e_H -e_K}
  \\
  &=n^{4v_G-3v_{H_0}}p^{4e_G-3e_{H_0}}\\\label{eq:max=max}
  &=\max_{ H \subset G \atop
    e_H \geq1}n^{4v_G-3v_H }p^{4e_G-3e_H },
\end{align}
as required. We proceed similarly with the sum $S_2$.
For $1\leq l< k\leq n$ we have 
\begin{align}
  \nonumber 
  \left\| g_l \star_{l}^{l} g_k\right\|^2_{\ell^2(\N)^{\otimes 2(k-l)}}
  &\approx\sum_{c\in\N^{k-l}}\(
  \sum_{b\in\N^l}\(\sum_{a\in\N^{e_G-l}}\mathbf1_{E_G}\(a,b\)\sum_{a'\in\N^{e_G-k}}\mathbf1_{E_G}\(a',b,c\)\)\)^2
    \\
\label{dssds0} 
    &\approx \sum_{A \subset K_n
      \atop
      e_A ={k-l}}\(
    \sum_{\substack{ A \subset B \subset K_n
        \\
        e_B =k}}\(
    \sum_{\substack{
        B \setminus A \subset G'' \subset K_n
        \atop
        G''\sim G}}
    1
    \sum_{\substack{
        H\subset G' \subset K_n
        \atop G'\sim G}}1
     \)\)^2
      \\
      \label{dssds1}
      &\lesssim
      \sum_{K\subset G
        \atop e_K=k-l}
      n^{v_K} n^{v_K}
      \(
      \sum_{\substack{K\subset H \subset G, 
          \ \! 
          H' \subset G
          \atop e_H =k, \ \! e_{H'}=l}}n^{v_H -v_K}\(n^{v_G-v_{H'}}n^{v_G-v_H}\)\)^2
        \\[10pt]
        \label{dssds2}
        &\lesssim \max_{\substack{K, H'\subset G\\e_K=k-l, \ \! e_{H'}=l}}n^{4v_G-2v_{H'}-v_K},
\end{align}
where $H'$ in \eqref{dssds1} stands for $B \setminus A$
in \eqref{dssds0}, whereas in \eqref{dssds2} the sum over $H'$ extends to all
$H'\subset G$ such that $e_{H'}=l$. It follows that 
\begin{align*}
  S_2 &\lesssim \sum_{1\leq l<k\leq e_G}p^{4e_G-k-l}\max_{\substack{K, H'\subset G
      \atop e_K=k-l, \ \! e_{H'}=l}}n^{4v_G-2v_{H'}-v_K}\\
  &=\sum_{1\leq l<k\leq e_G}\max_{\substack{K, H'\subset G
      \atop
      e_K =k-l, \ \! e_{H'}=l}}n^{4v_G-2v_{H'}-v_K}p^{4e_G-2v_{H'}-e_K}
  \\
  &\lesssim\max_{\substack{K', H'\subset G \atop
      e_{K'}, \ \!
      e_{H'}\geq1}}n^{4v_G-2v_{H'}-v_{K'}}p^{4e_G-2v_{H'}-e_{K'}}\\
&=n^{4v_G-3v_{H_0}}p^{4e_G-3e_{H_0}}\\
&=\max_{ H\subset G\atop e_H\geq1}n^{4v_G-3v_H}p^{4e_G-3e_H}, 
\end{align*}
    where $H_0$ is defined in \eqref{djlcn}. 
    Finally, we pass to estimates of $S_3$. For $1\leq l< k\leq n$ we have
\begin{eqnarray*}
  \lefteqn{
    \! \! \! \! \! \! \! \! \!
    \left\| g_k\star_{l}^{l} g_k\right\|^2_{\ell^2(\N)^{\otimes (k-l)}}
\approx\sum_{c,c'\in\N^{k-l}}\(
  \sum_{b\in\N^l}\(\sum_{a\in\N^{e_G-k}}\mathbf1_{E_G}\(a,b,c\)\)\(\sum_{a'\in\N^{e_G-k}}\mathbf1_{E_G}\(a',b,c'\)\)\)^2
  }
  \\
   & \approx  & \sum_{\substack{A ,A'\subset K_n
        \atop e_A =e_{A'}=k-l}}\(
    \sum_{\substack{ B \subset K_n
        \atop
        e_B =l, \ \! e_{ A \cap B }=e_{A' \cap B}=0}}\(\sum_{\substack{ A \cup B \subset G'\subset K_n
        \atop G'\sim G}}1\)
    \(\sum_{\substack{A'\cup B\subset G'' \subset K_n
        \atop G''\sim G}}1\)\)^2
  \\
      &\approx & \sum_{\substack{K,K',H\subset G\\e_K =e_{K'}=k-l, \ \! e_H =l\\e_{K\cap H}=e_{K'\cap H}=0}}
      \sum_{\substack{
          A, A' \subset K_n \\
          A\sim K \\ A'\sim K'}}\(
      \sum_{\substack{
          B \subset K_n \\
          B \sim H \\ A \cap B \sim K\cap H \\ A' \cap B \sim K'\cap H}}
      \(\sum_{\substack{ A \cup B \subset G'\subset K_n
          \atop G'\sim G}}1\)\(\sum_{\substack{A' \cup B \subset G''\subset K_n
          \atop G''\sim G}}1\)\)^2
      \\
      &\approx & \sum_{\substack{K,K',H\subset G\\e_K =e_{K'}=k-l, \ \! e_H =l\\e_{K\cap H}=e_{K'\cap H}=0}}
      \ \sum_{\substack{
          A, A' \subset K_n \\
          A\sim K \\ A'\sim K'}}\(
      \sum_{\substack{
          B \subset K_n \\
          B \sim H \\ A \cap B \sim K\cap H \\ A' \cap B \sim K'\cap H}}
      \(
      n^{v_G-v_{A \cup B}}
      \)\(
      n^{v_G-v_{A' \cup B}}
      \)\)^2. 
\end{eqnarray*}
Next, we note that given $A, A' \subset K_n$ it takes 
$$
v_B - v_{A\cap B} - v_{A'\cap B} + v_{A\cap A' \cap B}
=
v_H - v_{K\cap H} - v_{K'\cap H} + v_{A\cap A' \cap B}
$$ 
vertices to create any subgraph
$B \sim H$ such that $A \cap B \sim K\cap H$
and
$A' \cap B \sim K'\cap H$, with the bound 
$$
v_{A\cap A' \cap B}
\leq
\frac{1}{2} v_{A\cap A'}
+
\frac{1}{2} v_{A'\cap B}
=
\frac{1}{2} ( v_{A\cap A'} + v_{K'\cap H} ). 
$$
 Hence we have 
\begin{eqnarray*}
  \lefteqn{
      \left\| g_k\star_{l}^{l} g_k\right\|^2_{\ell^2(\N)^{\otimes (k-l)}}
  }
            \\
            &\lesssim &
            \! \! \! \! \! 
        \sum_{\substack{K,K',H \subset G\\e_K =e_{K'}=k-l, \ \! e_H =l\\e_{K\cap H}=e_{K'\cap H}=0}}
        \ \sum_{\substack{
            A, A' \subset K_n \\
            A \sim K \\ A'\sim K'
        }}
        \! \! \! \!
        \(
        n^{v_H-v_{K\cap H}-v_{K'\cap H}+
           ( v_{A\cap A'} + v_{K'\cap H} ) / 2
          }\(n^{v_G-v_{K\cup H}}\)\(n^{v_G-v_{K'\cup H}}\)\)^2. 
\end{eqnarray*}
 In order to estimate the above sum using powers of $n$, 
we need to consider the possible intersections
$A \cap A'$ for $A, A' \subset K_n$, as follows: 
\begin{eqnarray}
  \nonumber
  \lefteqn{
    \! \! \! \! \! \! \! \! \! \! \! \!
    \sum_{\substack{K,K',H\subset G\\e_K =e_{K'}=k-l, \ \! e_H =l\\e_{K\cap H}=e_{K'\cap H}=0}}
    \ \sum_{\substack{
            A, A' \subset K_n \\
            A \sim K \\ A'\sim K'
    }}
    n^{4v_G+2v_H-2v_{K\cap H}-v_{K'\cap H}
      +v_{A \cap A'} -2v_{K\cup H}-2v_{K'\cup H}}
  }
  \\
  \nonumber
  &\lesssim & \sum_{\substack{K,K',H\subset G\\e_K =e_{K'}=k-l, \ \! e_H =l\\e_{K \cap H}=e_{K' \cap H}=0}}\ \sum_{i=0}^{v_K}n^{v_K+v_{K'}-i} \ \! n^{4v_G+2v_H-2v_{K \cap H} - v_{K' \cap H}
    + i-2v_{K \cup H}-2v_{K' \cup H}}
  \\
  \label{aux3}
  &\lesssim & \sum_{\substack{K,K',H\subset G\\e_K =e_{K'}=k-l, \ \! e_H =l\\e_{K \cap H}=e_{K' \cap H}=0}}
  \ \! n^{v_K+v_{K'}+4v_G+2v_H-2v_{K \cap H}-v_{K' \cap H}-2v_{K \cup H}-2v_{K' \cup H}}. 
\end{eqnarray}
Furthermore we have  
\begin{eqnarray*}
  \lefteqn{
    \! \! \! \! \! \! \! \! \! \! \! \! \! \! \! \! \! \! \! \! \! \!
    \! \! \! \! \! \! \! \! \! \! \! \! \! \! \! \! \! \! \! \! \! \!
    v_K+v_{K'}+4v_G+2v_H-2v_{K \cap H}- v_{K' \cap H} - 2v_{K \cup H}-2v_{K' \cup H}
  }
  \\
&=&4v_G-v_K -v_H -v_{K'\cup H},
\end{eqnarray*}
so the sum \eqref{aux3} can be estimated as 
\begin{align*}
\sum_{\substack{K,K',H \subset G\\e_K =e_{K'}=k-l, \ \! e_H =l\\e_{K \cap H}=e_{K' \cap H}=0}}\ n^{4v_G-v_K -v_H -v_{K'\cup H}}\lesssim \max_{\substack{K,H,L\subset G\\e_K =k-l, \ \! e_H =l,\ \! e_L=k}} \ \! n^{4v_G-v_K -v_H -v_L },
\end{align*}
from which it follows
\begin{align*}
S_3&\lesssim \sum_{1\leq l<k\leq e_G}p^{4e_G-2k}
\max_{\substack{K,H,L\subset G\\e_K =k-l, \ \! e_H =l,\ e_L =k}}\ n^{4v_G-v_K -v_H -v_L }\\
&= \sum_{1\leq l<k\leq e_G}
\max_{\substack{K,H,L\subset G\\e_K =k-l, \ \! e_H =l,\ e_L =k}}\ n^{4v_G-v_K -v_H -v_L }p^{4e_G-e_K -e_H -e_L }\\
&\lesssim\max_{\substack{K,H,L \subset G\\e_K , e_H, e_L \geq1}}\ \! n^{4v_G-v_K -v_H-v_L }p^{4e_G-e_K -e_H-e_L }\\
&\leq n^{4v_G-3v_{H_0}}p^{4e_G-3e_{H_0}}\\
&= \max_{\substack{H\subset G\\e_H \geq1}}\ \! n^{4v_G-3v_H }p^{4e_G-3e_H },
\end{align*}
which ends the proof.
\end{Proof}
In the next corollary
 we note that Theorem~\ref{thm:main} simplifies 
 if we narrow our attention to $p_n$ depending of the complete
 graph size $n$ and close to $0$ or to $1$.
\begin{corollary}
\label{c0}
 Let $G$ be a graph without separated vertices. For $p_n<c<1$, $n\geq 1$, we have
\begin{equation}
\label{eq} 
  d_K \big(\widetilde{N}^G_n,\mathcal{N} \big) \lesssim
  \( \min_{\substack{ H\subset G\\e_H\geq1}}\big\{n^{v_H}p_n^{e_H}\big\}\)^{-1/2}.
\end{equation} 
 On the other hand, for  $p_n>c>0$, $n\geq 1$, it holds
 \begin{equation}
   \label{dk}
   d_K \big(\widetilde{N}^G_n,\mathcal{N} \big)\lesssim
 \frac{1}{n\sqrt{1-p_n}}.
\end{equation} 
\end{corollary}
 As a consequence of Corollary~\ref{c0}
 it follows that if 
 $$
 np_n^\beta \rightarrow \infty \ \mbox{ and }\ n^2(1-p_n)\rightarrow\infty,
 $$
 where $\beta :=\max \big\{e_H/v_H \ : \ H\subset G\big\}$,
 then we have the convergence of the
 renormalized subgraph count $\big(\widetilde{N}^G_n\big)_{n\geq 1}$
 to $\mathcal{N}$ in distribution as $n$ tends to infinity,
 which recovers the sufficient condition in \cite{rucinski}. 
 When $p\approx n^{-\alpha}$, $\alpha>0$,
 Corollary~\ref{c0} also shows that 
\begin{align}\label{eq:dFNwithalpha}
  d_K \big(\widetilde{N}^G_n,\mathcal{N} \big)
  &\lesssim
  \(\min_{\substack{ H\subset G\\e_H\geq1}}\big\{n^{v_H-\alpha e_H}\big\}\)^{-1/2}, 
\end{align}
 and in order for the above bound \eqref{eq:dFNwithalpha}
to tend to zero as $n$ goes to infinity, 
we should have 
\begin{equation}
  \label{alpha}
  \alpha<\min_{H\subset G}\frac{v_H}{e_H}=:\frac1\beta.
\end{equation} 
The next 
Corollary~\ref{cor1.3} of Theorem~\ref{thm:main} and \eqref{dk}
deals with cycle graphs with $r$ vertices, $r\geq3$.
When $G$ is a triangle it recovers the Kolmogorov bounds of \cite{roellin2}
as in Corollary~\ref{c} above.
\begin{corollary}
\label{cor1.3}
Let $G$ be a cycle graph with $r$ vertices, $r\geq3$, and
$c\in (0,1)$. We have 
\begin{align*}
  d_K \big(\widetilde{N}^G_n,\mathcal{N}\big)\lesssim
  \left\{\begin{array}{ll}
  \displaystyle
 \frac{1}{n\sqrt{1-p_n}} &\mbox{if } \ \displaystyle
  0<c< p_n, 
 \\
 \\
  \displaystyle
  \frac{1}{n\sqrt{p_n}} &
  \displaystyle
  \mbox{if } \ n^{-(r-2)/(r-1)} < p_n \leq c, 
  \\
  \\
  \displaystyle
 \frac{1}{(np_n)^{r/2}}&
  \displaystyle
  \mbox{if } \ 0 < p_n \leq n^{-(r-2)/(r-1)}.
  \end{array}\right.
\end{align*}
\end{corollary}
\begin{Proof}
 The smallest number of vertices of subgraphs $H$ of $G$ having $k$ edges,
 $k<r$, is realised for a linear subgraph having $k+1$ vertices,
 which yields 
 $$
  \min_{\substack{ H\subset G
      \atop
      1 \leq e_H < r}}
  \big\{n^{v_H} p_n^{e_H}\big\}
  =
  \min_{1\leq k < r}n^{k+1}p_n^k
  =n\min_{1\leq k < r} (n p_n )^k
  = \min (n^2p_n , ( n p_n )^{r-1}), 
  $$
  hence
  \begin{align*}
  \min_{\substack{ H\subset G\\e_H\geq1}} \big\{n^{v_H}p_n^{e_H}\big\}
  =\min\big\{n^2p_n,(np_n)^{r-1} , (np_n)^r
  \big\}
  =\left\{\begin{array}{ll}
  \displaystyle
 n^2p_n &\mbox{if } \  \displaystyle
  n^{-(r-2)/(r-1)} < p_n \leq c,
  \\
  \\
\displaystyle
(np_n)^r &\mbox{if } \ \displaystyle
0 < p_n \leq n^{-(r-2)/(r-1)},
  \end{array}\right.
\end{align*}
which concludes the proof by \eqref{eq} and \eqref{dk}. 
\end{Proof}
In case $p_n\approx n^{-\alpha}$
we should have $\alpha\in (0,1)$ by
  \eqref{alpha},
 Corollary~\ref{cor1.3} also shows that 
\begin{align*}
  d_K \big(\widetilde{N}^G_n,\mathcal{N} \big)\lesssim
  \left\{\begin{array}{ll}
  \displaystyle
  n^{-1+\alpha/2}
  \approx \frac{1}{n\sqrt{p_n}} &
  \displaystyle
  \mbox{if } \ 0 < \alpha\leq\frac{r-2}{r-1},
  \\
  \\
  \displaystyle
  n^{-r(1-\alpha)/2}
  \approx \frac{1}{(np_n)^{r/2}}&
  \displaystyle
  \mbox{if } \ \frac{r-2}{r-1}<\alpha<1.
  \end{array}\right.
\end{align*}
when $G$ is a cycle graph with $r$ vertices, $r\geq3$.
In the particular case $r=3$ where $G$ is a triangle, 
this improves on the Kolmogorov bounds in
Theorem~1.1 of \cite{reichenbachsAoP}. 
\\
 
In the case of complete graphs,
the next corollary also covers the case of
triangles. 
\begin{corollary}
\label{cor1.4}
 Let $G$ be a complete graph with $k$ vertices,
 $r\geq3$, and $c\in (0,1)$. We have 
\begin{align*}
  d_K \big(\widetilde{N}^G_n,\mathcal{N} \big)\lesssim
  \left\{\begin{array}{ll}
  \displaystyle
 \frac{1}{n\sqrt{1-p_n}} &\mbox{if } \ \displaystyle
  c< p_n <1, 
 \\
 \\
    \displaystyle
 \frac{1}{n\sqrt{p_n}} &
      \displaystyle
      \mbox{if } \
 n^{-2/(r+1)} < p_n \leq c, 
      \\
      \\
  \displaystyle
 \frac{1}{n^{r/2}p_n^{r(r-1)/4}}&  \displaystyle
 \mbox{if } \ 0 < p_n \leq n^{-2/(r+1)}. 
  \end{array}\right.
\end{align*}
\end{corollary}
\begin{Proof}
  The greatest number of edges of subgraphs of $G$ having $k$ vertices,
  $2\leq k\leq v_G$,
  is realised for a complete graph having $k\choose 2$ edges,
  which shows that 
  $$
  \min_{\substack{ H\subset G\\e_H\geq1}}\big\{n^{v_H}p_n^{e_H}\big\}=
  \min_{1\leq k\leq r}n^kp_n^{{k\choose 2}}.$$
On the other hand, from the equality 
$$
\frac{n^{k+1}p_n^{{k+1\choose 2}}}{n^k p_n^{{k\choose 2}}}=n p_n^k,
$$
we note that if the minimum was realised
with $k$ vertices where $1<k<r$,
we would have $n p_n^{k-1}\leq 1$ and
$np_n^k \geq1$, which would lead to 
$p_n \geq 1$, which is not possible. 
Therefore we have 
\begin{align*}
  \min_{\substack{ H\subset G\\e_H\geq1}}\big\{n^{v_H}p_n^{e_H}\big\}=\min\left\{
  n^2p_n,n^rp_n^{{r\choose 2}}\right\}=\left\{\begin{array}{ll}
  \displaystyle
  n^2p_n &\mbox{if } \ \displaystyle
  n^{-2/(r+1)} < p_n \leq c,
  \\
  \\ 
\displaystyle
  n^rp_n^{r(r-1)/2}&\mbox{if } \ \displaystyle
  0 < p_n \leq n^{-2/(r+1)},
  \end{array}\right.
\end{align*}
and we conclude the proof by \eqref{eq} and \eqref{dk}. 
\end{Proof}
When $p_n\approx n^{-\alpha}$ with $\alpha\in (0, 2/(r-1))$
by \eqref{alpha}, Corollary~\ref{cor1.4} shows that 
\begin{align*}
  \min_{\substack{ H\subset G\\e_H\geq1}}\{n^{v_H-\alpha e_H}\}=\min
  \big\{ n^{2-\alpha},n^{r-{r\choose 2}\alpha} \big\}=\left\{\begin{array}{ll}
\displaystyle
  n^{2-\alpha/2} &\mbox{if } \ \displaystyle
  0 < \alpha\leq\frac{2}{r+1},
  \\
  \\ 
 \displaystyle
  n^{r-r(r-1)\alpha/2}&\mbox{if } \ \displaystyle
  \frac{2}{r+1} \leq \alpha<\frac2{r-1},\end{array}\right.
\end{align*}
 hence by \eqref{eq} we find 
\begin{align*}
d_K \big(\widetilde{N}^G_n,\mathcal{N} \big)\lesssim
  \left\{\begin{array}{ll}
    \displaystyle
    n^{-1+\alpha/2} \approx \frac{1}{n\sqrt{p_n}} &
      \displaystyle
      \mbox{if } \ 0  < \alpha\leq\frac{2}{r+1},
      \\
      \\
  \displaystyle
  n^{-r/2+r(r-1)\alpha/4}
  \approx \frac{1}{n^{r/2}p_n^{r(r-1)/4}}&  \displaystyle
  \mbox{if } \ \frac{2}{r+1} \leq \alpha<\frac{2}{r-1}.\end{array}\right.
\end{align*}
Finally, the next corollary deals with
the important class of graphs which have a tree structure.
\begin{corollary}
  \label{cor1.5}
  Let $G$ be any tree (a connected graph without cycles) with $r$ edges,
  and $c\in (0,1)$.
  We have 
\begin{align*}
  d_K \big(\widetilde{N}^G_n,\mathcal{N} \big)\lesssim
  \left\{\begin{array}{ll}
    \displaystyle
 \frac{1}{n\sqrt{1-p_n}} &\mbox{if } \  \displaystyle
 c< p_n<1, 
 \\
 \\
    \displaystyle
   \frac{1}{n\sqrt{p_n}} &  \displaystyle
    \mbox{if } \ \frac{1}{n} < p_n \leq c,\\ \\
  \displaystyle
   \frac{1}{n^{(r+1)/2}p_n^{r/2}}&  \displaystyle
   \mbox{if } \ 0 < p_n \leq \frac{1}{n}.
  \end{array}\right.
\end{align*}
\end{corollary}
\begin{Proof}
  We have
  $$
  \min_{\substack{ H\subset G\\e_H\geq1}}\big\{n^{v_H} p_n^{e_H}\big\}=\min_{1\leq k\leq r}n^{k+1}p_n^k =n\min_{1\leq k\leq r} (n p_n )^k.$$
   The smallest number of vertices for a subgraph
  of a tree $G$ having $k$ edges, $k\leq r$, is
  realised for a subtree having $k+1$ vertices,
hence since $np_n$
  can be either less or greater than $1$, which gives 
$$\min_{\substack{ H\subset G\\e_H\geq1}}\big\{n^{v_H}p_n^{e_H}\big\}=n\min\big\{np_n,\(np_n\)^r\big\}=\left\{\begin{array}{ll}
\displaystyle
n^2 p_n & \mbox{if } \  \displaystyle
\frac{1}{n} < p_n \leq c,\\
\\
\displaystyle
n^{r+1}p_n^r &\mbox{if } \ \displaystyle
0 < p_n \leq \frac{1}{n},
  \end{array}\right.$$
as required, and we conclude by \eqref{eq} and \eqref{dk} .
\end{Proof}
In case $p_n\approx n^{-\alpha}$ with $\alpha \in ( 0, 1+1/r)$, we have
$\beta=\max \{e_H/v_H \ : \ H\subset G\} = r/(r+1)$ hence 
$$\min_{\substack{ H\subset G\\e_H\geq1}}\{n^{v_H-\alpha e_H}\}=n\min\{n^{1-\alpha},\(n^{1-\alpha}\)^r\}=\left\{\begin{array}{ll}
\displaystyle
n^{2-\alpha} & \mbox{if } \  \displaystyle
0 < \alpha\leq1,\\
\\
\displaystyle
n^{r+1-r\alpha} &\mbox{if } \ \displaystyle
1 \leq \alpha<1+\frac1r,\end{array}\right.$$
 which shows by \eqref{eq} that 
\begin{align*}
  d_K \big(\widetilde{N}^G_n,\mathcal{N} \big)\lesssim
  \left\{\begin{array}{ll}
    \displaystyle
    n^{-1+\alpha/2}
    \approx \frac{1}{n\sqrt{p_n}} &  \displaystyle
    \mbox{if } \ 0 < \alpha\leq1,
    \\ \\
  \displaystyle
  n^{-(r+1-r\alpha)/2}
  \approx \frac{1}{n^{(r+1)/2}p_n^{r/2}}&  \displaystyle
  \mbox{if } \ 1 \leq \alpha<1+\frac{1}{r}.
  \end{array}\right.
\end{align*}

\footnotesize

\def\cprime{$'$} \def\polhk#1{\setbox0=\hbox{#1}{\ooalign{\hidewidth
  \lower1.5ex\hbox{`}\hidewidth\crcr\unhbox0}}}
  \def\polhk#1{\setbox0=\hbox{#1}{\ooalign{\hidewidth
  \lower1.5ex\hbox{`}\hidewidth\crcr\unhbox0}}} \def\cprime{$'$}

\end{document}